\documentclass{article}
\usepackage{latexsym}
\usepackage{oldlfont}
\newtheorem{lemma}{Lemma}[section]
\newtheorem{theorem}[lemma]{Theorem}
\newtheorem{remark}[lemma]{Remark}

\newtheorem{coro}[lemma]{Corollary}
\newtheorem{definition}[lemma]{Definition}
\newcommand{\eps}{\varepsilon}
\newcommand{\Wsob}{\smash{{\stackrel{\circ}{W}}}_2^1(D)}

\newcommand{\om}{\omega}
\newcommand{\Om}{\Omega}

\renewcommand{\phi}{\varphi}

\input amssym.def
\input amssym
\title
{Determining functionals for random partial differential equations}
\author{Igor Chueshov\thanks{on leave from Department of Mechanics and
Mathematics,
Kharkov University, 310077 Kharkov, Ukraine}\\
Institute f\"{u}r Dynamische Systeme, FB3\\
Universit\"{a}t Bremen\\ D-28334 Bremen,
 Germany\\
\and
Jinqiao Duan\\
Department of Applied Mathematics\\
   Illinois Institute of Technology\\
   Chicago, IL 60616, USA \\
\and
Bj{\"o}rn Schmalfu{\ss}\\
Department of Applied Sciences\\
University of Technology and Applied Sciences\\
Geusaer Stra{\ss}e\\
D--06217 Merseburg,
Germany\\
}
\date{March 15, 2001}
\begin{document}
\maketitle
$\,$\\


\begin{abstract}

{\bf Nonlinear
    Diff. Eqns. Appl.} (NoDEA)  10(2003), 431-454.

Determining functionals are tools to describe the finite
dimensional long-term dynamics of infinite dimensional dynamical
systems. There also exist several applications to infinite
dimensional {\em random} dynamical systems. In these applications  the
convergence condition of the trajectories of an infinite
dimensional random dynamical system with respect to a finite set
of linear functionals is assumed to be either in mean or {\em exponential} with
respect to the convergence almost surely. In contrast
to these ideas we introduce a convergence concept which is based
on the convergence in probability. By this ansatz we get rid of
the assumption of exponential convergence. In addition, setting
the random terms
to
zero we obtain usual deterministic results.\\
We apply our results to the 2D Navier - Stokes equations
forced by a white noise.

\end{abstract}

\section{Introduction}
The question of the number of parameters that are necessary for
the description  of the long-term behaviour of solutions to
nonlinear partial differential equations was first discussed by
Foias and Prodi \cite{FP} and by Ladyzhenskaya \cite{La2} for the
deterministic 2D Navier-Stokes equations. They proved that the
asymptotic behaviour of the solutions is completely determined by
the dynamics of the first $N$ Fourier modes, if $N$ is
sufficiently large. After \cite{FP} and \cite{La2} similar results
were obtained for other parameters and other deterministic
equations and a general approach to the problem of the existence
of a finite number of determining parameters was developed (see
\cite{C3, C4, johtit93}  and the literature quoted therein).
\medskip\par
Assume that we have a dynamical system with the phase state $H$
and the evolution operator $S_t$. Roughly speaking the general
problem  on the existence of finite sets of determining parameters
(functionals) can be stated (cf. \cite{C3, C4}) as follows: find
the conditions on a finite set $\{ l_j : j= 1,...,N\}$  of
functionals on $H$ which guarantees that the convergence (in
certain sense)
\[
\max_j|l_j(S_tu_1-S_tu_2)|\to 0\quad\mbox{when}\quad t\to+\infty,\; u_1, u_2\in H
\]
implies that $S_tu_1-S_tu_2\to 0$ in some topology of the phase space $H$.
Besides from applied point of view it is also important to find bounds for
the number $N$ of determining functionals (in the sense above) and to describe families
of functionals with minimal $N$.
\medskip\par
The deterministic theory of determining functionals was developed
by many authors (see, e.g.   \cite{C3, C4}  and the references
therein). Similar   problems  for stochastic systems were also
discussed in \cite{C2, BF, FL, C5}. In papers   \cite{C2, BF, FL}
$\om$-wise approach to construction of determining functionals
were developed. However in these papers it was assumed that either
(see \cite{C2, BF, FL}) the nonlinear term in the equation is a
globally Lipschitz mapping or (see \cite{BF, FL}) one of initial
data belongs to the random attractor. The mode of convergence for
functionals and trajectories is the convergence
 almost surely
in these papers. Moreover in the papers \cite{BF, FL}  the authors assume that the
functionals of the difference of two solutions go to zero exponentially fast.
Then they prove that some norm of the difference of these solutions tends to zero
with exponential speed which is less than the speed of convergence of the functionals.
On the other hand the approach presented in \cite{C5} does not assume  these conditions, and
it relies on some estimates exponential moments of solutions and deals with
convergence  in mean. The speed of convergence to zero of the functionals and
of the norms are the same as in \cite{C5}.
\medskip\par
In contrast to  \cite{ C2, BF, FL, C5} we consider determinig functionals
with respect to the convergence in probability. Using such determining functionals
we can avoid the assumption that the images (under linear functionals)
of the trajectories converge
 exponentially fast. In particular, our approach
recovers the deterministic results when we remove the stochastic terms.
Another advantage is that we do not have   to assume that  one trajectory
must be contained in the random attractor. Finally, we mention that
the convergence in probability is quite natural for RDS,
see \cite{CraFla94, ArnSchm99}.
\medskip\par
In Section~2 we consider an abstract setting of random dissipative
systems and prove two
 existence theorem of finite sets of determining functionals
in the sense of the definition given below for arbitrary initial data.
These theorems show  two different approaches to the construction of
 determining functionals.
In Section~3 we apply the results of Section~2 to
2D Navier - Stokes equations subject  to additive white noise.
We prove the existence of
determining functionals for this problem {\em without}
the assuming that one of solutions belongs to the attractor.

\section{The existence of determining functionals}
We consider a {\em random dynamical system}  (REDS) which consists
of two components. The first component is a {\em metric dynamical
system} $(\Omega,{\cal F},\Bbb{P},\theta)$ as a model for a
noise, where $(\Omega,{\cal F},\Bbb{P})$ is a probability
space and $\theta$ is a $\cal{F}\otimes
{\cal B}(\Bbb{R}),{\cal F}$ measurable flow: we have
\[
\theta_0={\rm id},\qquad
\theta_{t+\tau}=\theta_t\circ\theta_\tau=:\theta_t\theta_\tau
\]
for $t,\,\tau\in\Bbb{R}$. The measure $\Bbb{P}$ is supposed
to be ergodic with respect to $\theta$. The second  component of a
random dynamical system is a
${\cal B}(\Bbb{R}^+)\otimes{\cal F}\otimes{\cal B}(H),{\cal B}(H)$-measurable
mapping $\phi$ satisfying the {\em cocycle} property
\[
\phi(t+\tau,\omega,x)=\phi(t,\theta_\tau\omega,\phi(\tau,\omega,x)),\qquad \phi(0,\omega,x)=x,
\]
where the phase space $H$ is a separable metric space and $x$ is chosen arbitrarily in $H$. We will denote this RDS by symbol $\phi$.\\
A standard model for such a noise $\theta$ is the twosided {\em
Brownian motion}: Let $U$ be a separable Hilbert space. We
consider the probability space
\[
(C_0(\Bbb{R},U),{\cal B}(C_0(\Bbb{R},U)),\Bbb{P})
\]
where $C_0(\Bbb{R},U)$ is the Fr\'echet space of continuous
functions on $\Bbb{R}$ which are zero at zero and
${\cal B}(C_0(\Bbb{R},U))$ is the corresponding Borel
$\sigma$-algebra. Suppose that we have a covariance operator $Q$
on $U$. Then $\Bbb{P}$ denotes the {\em Wiener measure} with
respect to $Q$. Note that $\Bbb{P}$ is ergodic with respect to
the flow
\[
\theta_t\omega=\omega(\cdot+t)-\omega(t),\qquad\mbox{for
}\omega\in C_0(\Bbb{R},U).
\]
For  detailed presentation of random dynamical systems we
refer to the monograph by L. Arnold \cite{Arn98}.
\medskip\par\noindent
On a $V\subset H\subset V^\prime$  rigged Hilbert  space with compact
embedding $V\subset H$ and duality mapping $\langle\cdot,\cdot\rangle$
we investigate RDS $\phi$ generated by the
evolution equation
\begin{equation}\label{eq-11}
\frac{du}{d\,t}+A\,u=F(u,\theta_t\omega),\quad u(0)=x,
\end{equation}
over  some   metric dynamical system
$(\Omega,{\cal F},\Bbb{P},\theta)$. Here $A$ is a positive
self-adjoint operator in $H$ such that $D(A^{1/2})=V$, where
$D(B)$ denotes the domain of the operator $B$. We suppose that $V$
is equipped with the norm $\|\cdot\|_V =\| A^{1/2}\cdot\|_H$. We
also assume that the nonlinear mapping $F$ from $V\times \Omega$
into $H$ is such that $F(u,\omega)$ is measurable for any fixed
$u\in V$ and subordinate (in the sense of (\ref{eqb1})) to the
operator $A$. We suppose that the solution $u(t,\omega)$ of the
problem (\ref{eq-11}) is unique and depends  measurably on
$(t,\omega,x)$. Then the operator
\[
(t,\omega,x)\to u(t,\omega),\qquad u(0,\omega)=x
\]
 defines a
random dynamical system (cocycle) $\phi$. In addition, this random dynamical
system is supposed to be continuous which means that
\[
x\to \phi(t,\omega,x)
\]
is continuous for any $(t,\omega)$. The trajectories of this random dynamical system
has to be contained in $L_{2,loc}(0,\infty;V)\cap C([0,\infty);H)$.\\

In the following we assume that $\phi$ is {\em dissipative}. It means that there
exists a {\it compact } random set $B\subset V $ which
is forward invariant:
\[
\phi(t,\omega,B(\omega))\subset B(\theta_t\omega)
, \; t>0,
\]
and which is absorbing: for any $\eps>0$ and for any random variable $x(\omega)\in H$
there exists a $t_{\eps,x}>0$ such that if $t\ge t_{\eps,x}$
\[
\phi(t,\omega,x(\omega))\in B(\theta_t\omega)
\]
with probability $1-\eps$. Note that $B$ is absorbing with probability one, due to
the forward invariance.
\\
A random variable $x\ge 0$ is called tempered if
\[
\lim_{t\to\pm\infty}\frac{\log^+x(\theta_t\omega)}{|t|}=0
\quad\mbox{a.s.}
\]
Note that the only alternative to this property is that
\[
\limsup_{t\to\pm\infty}\frac{\log^+x(\theta_t\omega)}{|t|}=\infty\quad\mbox{a.s.},
\]
see Arnold \cite{Arn98}, page 164 f.
We also assume that $B$ is tempered which means that the mapping
\[
\omega\to\sup_{x\in B(\omega)}\|x\|_H
\]
is tempered.

We now give our basic definition:
\begin{definition}
A set ${\cal L}=\{l_j,\;j=1,\cdots, k\}$ of linear continuous
and linearly independent functionals on $V$ is called
asymptotically determining in probability
if
\[
(\Bbb{P})\lim_{t\to\infty}\int_t^{t+1}
\max_j|l_j(\phi(\tau,\omega,x_1)-\phi(\tau,\omega,x_2))|^2d\tau\to 0
\]
for two initial conditions $x_1,\,x_2\in H$ implies
\[
(\Bbb{P})\lim_{t\to\infty}\|\phi(t,\omega,x_1)-\phi(t,\omega,x_2)\|_H\to
0.
\]
\end{definition}
As in \cite{C3, C4} we use the concept of the {\em completeness
defect} for a description of sets of determining functionals.
Assume that $X$ and $Y$ are Banach spaces and $X$ continuously and
densely embedded into $Y$.
Let ${\cal L} =  \{ l_j : j= 1,...,k\}$ be a finite set of
linearly independent continuous functionals on
$X$.  We define the completeness defect
$\eps_{\cal L}(X, Y)\equiv \eps_{\cal L}$
of the set ${\cal L}$ with respect to the pair of the spaces
$X$ and $Y$ by the formula
$$
 \eps_{\cal L} = \sup\{ \Vert w \Vert_Y \; : w\in X, \;
 l_j (w)= 0,\; l_j\in {\cal L}, \; \Vert w \Vert_X \le 1 \} .
$$
The value $\epsilon_{\cal L}$ is proved to be very useful for
characterization of sets of determining functionals (see, e.g.,
\cite{C3, C4} and the references therein). One can show that the
completeness defect $\eps_{\cal L}(X, Y)$ is the best possible
global error of approximation in $Y$ of elements $u\in X$ by
elements of the form $u_{\cal L}= \sum_{j= 1}^k l_j(u)\phi_j$,
where  $\{\phi_j : j= 1,\ldots , k\}$ is an arbitrary set in $X$.
The smallness of $\eps_{\cal L}(X, Y)$ is the main condition
(see the results presented below) that guarantee the property of a
set of functionals to be asymptotically determining. The so-called
modes, nodes and local volume averages (the description of these
functionals  can be found  in \cite{C4}, for instance) are the
main examples of sets of functionals with a small completeness
defect. For further discussions and for other  properties of the
completeness defect we refer to \cite{C3, C4}. Here we only point
out the following estimate
\begin{equation}\label{eqd}
 \Vert u\Vert_Y \le
\eps_{\cal L} \cdot \Vert u \Vert_X
+  C_{\cal L} \cdot\max_{j= 1,\ldots, k} \vert l_j (u) \vert,\quad u\in X,
\end{equation}
where $C_{\cal L}>0$ is a constant depending on ${\cal L}$.
\medskip

We are now   in a position to prove the first main theorem
for systems introduced in (\ref{eq-11}).
To do this we will use  the completeness defect
$\eps_{\cal L}\equiv\eps_{\cal L}(X, Y)$ with $H=Y,\,V=X$.
\begin{theorem}\label{t1}
Let ${\cal L} =\{ l_j : j= 1,...,k\}$ be  a set of linear continuous
and linearly independent functionals on $V$.
We assume that
we have a forward absorbing and forward invariant set $B$ in $V$ such that
$\sup_{x\in B(\omega)}\|x\|_V^2$
is bounded by a tempered random variable
and $t\to\sup_{x\in B(\theta_t\omega)}\|x\|_V^2$
is locally integrable.
Suppose there exist a constant $c>0$ and a measurable function $l\ge0$ such that for
$x_1(\omega),\,x_2(\omega)\in B(\omega)$ we have
\begin{equation}\label{eqb1}
\langle - A(x_1-x_2)+
F(x_1,\omega)-F(x_2,\omega),x_1-x_2\rangle
\end{equation}
\[
\le
-c\|x_1-x_2\|_V^2+ l(x_1,x_2,\omega)\|x_1-x_2\|_H^2.
\]
Assume that
\begin{equation}\label{eqb2}
\frac{1}{m}\Bbb{E}\left\{\sup_{x_1,x_2\in B(\omega)}\int_0^m
l(\phi(t,\omega,x_1),\phi(t,\omega,x_2),\theta_t\omega)dt\right\}<
c\eps_{{\cal L}}^{-2}
\end{equation}
for some $m>0$.
Then ${\cal L}$ is a set of asymptotically determining functionals in probability for RDS $(\theta,\phi)$.
\end{theorem}

{\it Proof.}
1) Without loss of generality, we only consider the case $m=1$. That is, we assume that
 (\ref{eqb2}) is fulfilled for
$m=1$.
Since we intend to prove  convergence in probability
we can suppose that the  random variables $x_1(\omega),\,x_2(\omega)$ are
contained in  $B(\omega)$.
Such random variables exist because $B$ is a random set, see Caistaing and Valadier
\cite{CasVal77} Chapter III. Indeed, $B$ is {\em forward} absorbing such that
$\phi(t,\omega,x_i(\omega))\in B(\theta_t\omega)$ with probability
$1-\eps$ for any $\eps >0$ if $t$ is sufficiently large.\\

Let $w(t,\omega)$ be defined by
$\phi(t,\omega,x_1(\omega))-\phi(t,\omega,x_2(\omega))$. Since
$\|\cdot\|_V =\| A^{1/2}\cdot\|_H$, we obtain by (\ref{eqb1}):
\[
\frac{d\|w\|_H^2}{dt}+2c\|w\|_V^2 \le
2l(\phi(t,\omega,x_1(\omega)),\phi(t,\omega,x_2(\omega)),\theta_t\omega)\|w\|_H^2\, .
\]
We have by (\ref{eqd})
\[
\|w\|_V^2\ge(1+\delta)^{-1}\eps_{{\cal L}}^{-2}\|w\|_H^2
-C_{\delta,{\cal L}}\max_{j=1,\cdots,k}|l_j(w)|^2.
\]
for any $\delta>0$ with appropriate positive constant
$C_{\delta,{\cal L}}$.
This allows us
to write the following inequality:
\begin{equation}\label{eqe}
\|w(t)\|_H^2\le\|w(0)\|_H^2e^{\int_0^t q(s,\omega)ds}\
+C_{\delta,{\cal L}}\cdot\int_0^t
e^{\int_s^tq(\tau,\omega)d\tau}\eta_{{\cal L}}(s,\omega)ds\ ,
\end{equation}
where $\eta_{{\cal L}}(s,\omega)=\max_j|l_j(w(s,\omega))|^2$ and
\[
q(t,\omega)=
2l(\phi(t,\omega,x_1(\omega)),\phi(t,\omega,x_2(\omega)),\theta_t\omega)
-2c(1+\delta)^{-1}\eps_{{\cal L}}^{-2}.
\]
Let
\[
Q(\omega)=\sup_{x_1,x_2\in B(\omega)}\int_0^1
2l(\phi(t,\omega,x_1),\phi(t,\omega,x_2),\theta_t\omega)dt-
2\,c(1+\delta)^{-1}\eps_{{\cal L}}^{-2}
\]
with $\delta>0$ chosen such that $\Bbb{E}Q<0$.
This is possible because of (\ref{eqb2}).\\

2) Since $B$ is forward invariant and
$q(t,\om)\ge -2c(1+\delta)^{-1}\eps_{{\cal L}}^{-2}$,
the first term on the right hand side
of (\ref{eqe}) can be estimated by
\[
\|w(0)\|_H^2e^{\sum_{j=0}^{[t]}Q(\theta_j\omega)}e^{2c(1+\delta)^{-1}\eps_{{\cal L}}^{-2}}.
\]
Since $\Bbb{E}Q<0$, by the ergodic theorem we have for
$t\to\infty$
\begin{equation}\label{eqg}
\sum_{j=0}^{[t]}Q(\theta_j\omega)\sim ([t]+1)\Bbb{E}Q\to-\infty
\end{equation}
which shows the convergence assertion for the first term. \\
We now investigate the second term in (\ref{eqe}). Since $B$ is
forward invariant, this  term can be estimated by
\begin{eqnarray*}
&&
C_{\delta,{\cal L}}\cdot\int_0^{[t]+1}
e^{\int_{[s]}^{[t]+1}q(\tau,\omega)d\tau}\eta_{{\cal L}}(s,\omega)ds
\,e^{4c(1+\delta)^{-1}\eps_{{\cal L}}^{-2}}\\
&&
\qquad\le C_{\delta,{\cal L}}\cdot
e^{4c(1+\delta)^{-1}\eps_{{\cal L}}^{-2}}
\sum_{j=0}^{[t]}e^{\sum_{j^\prime=j}^{[t]}Q(\theta_{j^\prime}\omega)}\int_0^1\eta_{{\cal L}}(s+j,\omega)ds\ .
\end{eqnarray*}
We use here  that $l(x_1, x_2,\omega)$ is a nonnegative function.
Thus we have to prove that
\begin{equation}\label{eqg1}
(\Bbb{P})\,\lim_{k\to\infty}
\sum_{j=0}^{k}\left(e^{\sum_{j^\prime=j}^{k}Q(\theta_{j^\prime}\omega)}\int_0^1\eta_{{\cal L}}(s+j,\omega)ds\right)=0.
\end{equation}
We now replace $\omega$ by $\theta_{-k}\omega$ in the relation under the
limit sign. It gives
\[
\sum_{j=-k}^{0}e^{\sum_{j^\prime=j}^{0}Q(\theta_{j^\prime}\omega)}\int_0^1\eta_{{\cal L}}(s+k+j,\theta_{-k}\omega)ds ,
\]
which is equal to
\[
\sum_{j=-\infty}^{0}e^{\sum_{j^\prime=j}^{0}Q(\theta_{j^\prime}\omega)}\chi_k(j)\int_0^1\eta_{{\cal L}}(s+k+j,\theta_{-k}\omega)ds
\]
where $\chi_k(j)=1$ if $j\ge-k$ and 0 otherwise.\\
3) Since we can assume that $x_i(\omega)\in B(\omega)$ there exists a
tempered random
variable $b$ such that
$\eta_{{\cal L}}(s,\omega)\le b(\theta_s\omega)$ where $s\to b(\theta_s\omega)$ is
locally integrable.
Since $s\to b(\theta_s\omega)$ is tempered
\[
j\to \chi_k(j)\int_0^1\eta_{{\cal L}}(s+k+j,\theta_{-k}\omega)\le \int_0^1b(\theta_{s+j}\omega)ds
\]
has a subexponential growth for any $k\ge 0$. We consider the
finite measure
$\mu^\omega(j)=e^{\frac{1}{2}\sum_{j^\prime=j}^0Q(\theta_{j^\prime}\omega)}\delta_{j}$
on $\Bbb{Z}^-$ where $\delta_j$ are Dirac measures on $j$. Set
\[
f(k,j,\omega):=\chi_k(j)e^{\frac{1}{2}\sum_{j^\prime=j}^0Q(\theta_{j^\prime}\omega)}\int_0^1\eta_{{\cal L}}(s+k+j,\theta_{-k}\omega)ds .
\]
Since $j\to\int_0^1\eta_{{\cal L}}(s+k+j,\theta_{-k}\omega)ds$ has a subexponential growth
and
$$
j\to e^{\frac{1}{2}\sum_{j^\prime=j}^0Q(\theta_{j^\prime}\omega)}
$$ goes to zero exponentially fast (see (\ref{eqg})),
there exists a constant
$n$   depending only on $\omega$ such that
\begin{equation}\label{eqh}
f(k,j,\omega)\le n(\omega)\qquad \mbox{for any}\quad -j,\,k\in
\Bbb{Z}^+.
\end{equation}
The term $\int_0^1\eta_{{\cal L}}(s+k+j,\omega)ds$ tends to
zero in probability for $k\to\infty$ and fixed $j$. Hence, also
$\int_0^1\eta_{{\cal L}}(s+k+j,\theta_{-k}\omega) ds$ tends to
zero in probability for $k\to\infty$ and fixed $j$. Let
$\lambda(j)=e^{\frac{1}{4}\Bbb{E}Qj}\delta_j,\,j\in
\Bbb{Z}^-$ be a finite measure on $\Bbb{Z}^-$ and
$\Xi_N(\omega)$ be the indicator function of the set
\[
\{ \omega:
e^{\frac{1}{2}\sum_{j^\prime=j}^0Q(\theta_{j^\prime}\omega)}\le
Ne^{\frac{1}{4}\Bbb{E}Qj} \;\mbox{for } j\in\Bbb{Z}^-\} ,
\]
where $\Xi_N$ tends increasingly to one for $N\to\infty$. We set
$\mu_N^\omega=\Xi_N(\omega)\mu^\omega$. For the asserted
convergence we consider
\begin{eqnarray*}
&&
\Bbb{P}(\int f(k,j,\omega)d\mu^\omega(j)>3\delta)\\
&& =\Bbb{P}(\int  (f(k,j,\omega)\wedge N) d\mu^\omega_N(j) +
\int f(k,j,\omega)-(f(k,j,\omega)\wedge N)d\mu^\omega_N(j)\\
&&\qquad
+\int f(k,j,\omega)d(\mu^\omega-\mu^\omega_N)(j)
>3\delta)\\
&& \le \Bbb{P}(\int (f(k,j,\omega)\wedge
N)d\mu^\omega_N(j)>\delta)\\
&&\qquad+ \Bbb{P}(\int f(k,j,\omega)-(f(k,j,\omega)\wedge
N)d\mu^\omega(j)
>\delta)\\
&& \qquad +\Bbb{P}(\int
f(k,j,\omega)d(\mu^\omega-\mu^\omega_N)(j)>\delta).
\end{eqnarray*}
Note that by (\ref{eqh}) the second term on the right hand side is
less than $\eps$ uniformly in $k$ if $N$ is sufficiently large
uniformly in $k$. The integral in the third term can be estimated
by
\[
\int f(k,j,\omega)d(\mu^\omega-\mu_N^\omega)(j)\le
n(\omega)(\mu^\omega-\mu^\omega_N)(\Bbb{Z}^-)=n(\om)(1-\Xi_N(\om))\mu^\om(\Bbb{Z}^-).
\]
Therefore this third term is less than $\eps$ for large $N$.
To see that the first term tends to zero in probability
for $k\to \infty$ and any $N$, we note that by the definition of the
metric of the convergence  in probability
\begin{eqnarray*}
&&\Bbb{E}\frac{\int  (f(k,j,\omega)\wedge N)
d\mu_N^\omega(j)}{1+\int  (f(k,j,\omega)\wedge N)
d\mu_N^\omega(j)}\le \Bbb{E}\int (f(k,j,\omega)\wedge
N)d\mu_N^\omega(j)
\\
&&\qquad = \int\Bbb{E} (f(k,j,\omega)\wedge N)d\mu_N^\omega(j)
\le N(N+1)\int\frac{\Bbb{E}(f(k,j,\omega)\wedge
N)}{1+\Bbb{E}(f(k,j,\omega)\wedge N)}d\lambda(j),
\end{eqnarray*}
where the right hand side tends to zero for any $N\ge 0$ by Lebesgue's theorem.
Thus   the asserted convergence (\ref{eqg1}) follows.
\hfill $\Box$\\

Now we present another version of the theorem on the existence of finite
number of determining functionals which can be easily applied to the
random squeezing property introduced by Flandoli and Langa \cite{FL}.

\begin{theorem}\label{t1a}
Let  $\phi$ be RDS whose phase space is a Banach space $H$ with
the norm $\Vert\cdot\Vert$. Suppose  that  this RDS is dissipative
in $H$ with a forward invariant absorbing random set $B(\om)$ such
that the random variable $\rho(\om)=\sup_{x\in B(\om)}\| x\|$ is
tempered and $\rho(\theta_t\om)\in L^p_{loc}({\Bbb R})$ for
some $p\ge 1$ and all $\om\in\Omega$. Assume that for each
$\om\in\Omega$ RDS $\phi$ possesses the following properties:
\begin{equation}\label{eqh1}
\Vert\phi(t,\om,x_1)-\phi(t,\om,x_2)\Vert\le M(\om) \Vert
x_1-x_2\Vert
\end{equation}
for all $t\in [0,1],\; x_1, x_2\in B(\om)$ and
\begin{eqnarray}\label{eqh2}
\Vert\phi(1,\om,x_1)-\phi(1,\om,x_2)\Vert&\le& {\cal
N}(\phi(1,\om,x_1)-\phi(1,\om,x_2))\nonumber\\
&+& e^{\int_0^1 r(\theta_\tau\om)d\tau}\cdot \Vert x_1-x_2\Vert
\end{eqnarray}
for all  $x_1, x_2\in B(\om)$. Here $M(\om)$ is a tempered and
finite almost surely random variable, ${\cal N}(\cdot)$ is a
positive continuous scalar function on $H$ such that  ${\cal
N}(x)\le C\cdot (1+\| x\|^p)$ and $r(\om)$ is a random variable
with finite expectation
 such that $\Bbb{E}r<0$. Then the condition
\begin{equation}\label{eqh3}
(\Bbb{P}) \lim_{n\to+\infty} {\cal
N}(\phi(n,\om,x_1)-\phi(n,\om,x_2))=0
\end{equation}
for some $x_1, x_2\in H$ implies that
\begin{equation}\label{eqh4}
(\Bbb{P}) \lim_{t\to+\infty} \Vert
\phi(t,\om,x_1)-\phi(t,\om,x_2)\Vert=0.
\end{equation}
\end{theorem}

{\it Proof.$\;$} As above we can assume that $x_i(\om)\in B(\om)$.
Using the cocycle property
$\phi(m,\om)=\phi(1,\theta_{m-1}\om,\phi(m-1,\om))$ and
relation (\ref{eqh2}) we obtain that
\[
d_m(\om)\le
{\cal N}(m,\om)+
e^{\int_{m-1}^m r(\theta_\tau\om)d\tau}\cdot d_{m-1}(\om),
\]
where
\[
{\cal N}(n,\om)=
{\cal N}(\phi(n,\om,x_1(\om))-\phi(n,\om,x_2(\om)))
\]
and
\[
d_t(\om)=\Vert \phi(t,\om,x_1(\om))-\phi(t,\om,x_2(\om))\Vert.
\]
After iterations we find that
\[
d_m(\om)\le d_{0}(\om)e^{\int_{0}^{m} r(\theta_\tau\om)d\tau}
+
\sum_{j=0}^{m-1} {\cal N}(m-j,\om)
e^{\int_{m-j}^m r(\theta_\tau\om)d\tau}.
\]
Applying now the same arguments as in the proof of Theorem
\ref{t1} we find that (\ref{eqh3}) implies that $(\Bbb{P})\,
\lim_{m\to+\infty}d_m(\om)=0$.
 From (\ref{eqh1}) we have that
$d_t(\om)\le M(\theta_{[t]}\om)d_{[t]}(\om)$. Thus we should prove
that $(\Bbb{P})\,
\lim_{n\to+\infty}M(\theta_{n}\om)d_{n}(\om)=0$. It follows from
$(\Bbb{P})\, \lim_{n\to+\infty}M(\om)d_{n}(\theta_{-n}\om)=0$.
The last relation follows from the convergence $(\Bbb{P})\,
\lim_{m\to+\infty}d_m(\theta_{-m}\om)=0$ and the properties of
$M(\om)$.
\hfill $\Box$\\

Now following   Flandoli and Langa \cite{FL} we introduce the
concept of random squeezing property.
\begin{definition}
Let  $\phi$ be RDS whose phase space is a separable Hilbert space
$H$. We say that RDS  $(\theta, \phi)$ satisfies a {\em random
squeezing property} (RSP) on the random set $B(\om)$ if there
exist a finite-dimensional projector $P$ and a random variable
$r(\om)$ with finite expectation such that $\Bbb{E}r<0$ and for
almost all $\om\in\Om$ we have either
\[
\Vert(I-P)\phi(1,\om,x_1)-\phi(1,\om,x_2)\Vert\le \Vert
P\phi(1,\om,x_1)-\phi(1,\om,x_2)\Vert
\]
or
\[
\Vert\phi(1,\om,x_1)-\phi(1,\om,x_2)\Vert\le
e^{\int_0^1 r(\theta_\tau\om)d\tau}\cdot
\Vert x_1-x_2\Vert
\]
for all  $x_1, x_2\in B(\om)$.
\end{definition}
In deterministic case a similar property is well-known for
dissipative systems with finite-dimensional long-time behaviour
(see, e.g. \cite{Tem97} and the references therein). Flandoli and
Langa \cite{FL} have proved random squeezing property  for a class of stochastic
reaction-diffusion equations and for stochastic $2D$
Navier~-~Stokes equations with periodic boundary condition.
\medskip\par\noindent
Now we are in position to state corollaries from Theorem~\ref{t1a}.
\begin{coro}\label{coro01}
Assume that  RDS $\phi$ with Hilbert phase space $H$
is dissipative with a forward
invariant absorbing random set $B(\om)$ satisfying the hypotheses of Theorem~\ref{t1a}.
Suppose that $\phi$
 possesses  property (\ref{eqh1}) and satisfies RSP
with  an orthogonal projector $P$. Then the property
\[
(\Bbb{P}) \lim_{n\to+\infty}
\left\{(\phi(n,\om,x_1),e_i)_H-(\phi(n,\om,x_2),e_i)_H\right\}=0,
\quad i=1,2,\ldots, d,
\]
for some $x_1, x_2\in H$ implies  (\ref{eqh4}).
Here $\{ e_i : i=1,\ldots,d\}$ is a basis in the subspace $PH$.
\end{coro}
{\it Proof.$\;$}
It is clear that RSP implies (\ref{eqh2}) with
${\cal N}(u) =2\Vert Pu\Vert$. Thus we can apply Theorem~\ref{t1a}.
\hfill $\Box$\\

This result on determining modes extends in some sense the result by
Flandoli and Langa \cite{FL} for the case $k=0$.
\begin{coro}\label{coro02}
Assume that  RDS $\phi$
satisfies the hypotheses of Corollary~\ref{coro01}.
Suppose that there exists a Banach space $W$ such that $H$ continuously and
densely embedded into $W$
and the projector $P$ can be extended to continuous operator from $W$
into $H$ such that  $\Vert Pu\Vert_H\le a_0 \Vert u\Vert_W$
with a positive constant $a_0$.
Let ${\cal L} =  \{ l_j : j= 1,...,k\}$ be a  set of
linearly independent continuous functionals on
$H$ with the completeness defect
$\eps_{\cal L}(H, W)$  with respect to the pair of the spaces
$H$ and $W$. If
\[
2 a_0\eps_{\cal L}(H, W)<1\quad\mbox{and}\quad
\Bbb{E}r+\log\frac{1}{1-2 a_0\eps_{\cal L}(H, W)}<0,
\]
then the property
\[
(\Bbb{P}) \lim_{n\to+\infty} \left\{
l_i(\phi(n,\om,x_1))-l_i(\phi(n,\om,x_2))\right\}=0 \quad
i=1,2,\ldots, k,
\]
for some $x_1, x_2\in H$ implies  (\ref{eqh4}).
\end{coro}
{\it Proof.$\;$}
As above  RSP implies (\ref{eqh2}) with
${\cal N}(u) =2\Vert Pu\Vert_H$. However using (\ref{eqd}) with $X=H$
and $Y=W$ we have
\[
2\Vert Pu\Vert_H\le 2a_0\Vert u\Vert_W \le
2 a_0\eps_{\cal L}(H, W)\cdot\Vert u\Vert_H +
C_{\cal L}\bar{\eta}(u),
\]
where
$\bar{\eta}_{\cal L}(u)=
\max\{ \vert l_j (u) \vert\, : j= 1,\ldots, k\}$.
Therefore from (\ref{eqh2}) we have
\begin{eqnarray*}
\Vert\phi(1,\om,x_1)&-&\phi(1,\om,x_2)\Vert_H\le \frac{C_{\cal
L}}{ 1-2 a_0\eps_{\cal L}(H, W)}\cdot
\bar{\eta}(\phi(1,\om,x_1)-\phi(1,\om,x_2))\\
&+& \exp\left\{ \int_0^1 r(\theta_\tau\om)d\tau+ \log\frac{1}{1-2
a_0\eps_{\cal L}(H, W)}\right\} \cdot \Vert x_1-x_2\Vert_H.
\end{eqnarray*}
Thus we can apply Theorem~\ref{t1a}.
\hfill $\Box$\\

\begin{remark}
{\rm
The space $W$ with the properties listed in Corollary~\ref{coro02}
can be easily constructed in the following situation. Assume that
$A$ is a positive self-adjoint operator in $H$ with compact resolvent.
Let $0<\lambda_1\le \lambda_2\le \ldots$ be the corresponding eigenvalues.
If $P$ is orthogonal projector on the first $k$ eigenvectors of $A$,
then we have $\Vert Pu\Vert_H\le \lambda_{k+1}^s \Vert A^{-s/2} u\Vert_H$
for any $s>0$.
Thus we can choose $W$ as a completion  of $H$  with respect to the norm
$\Vert A^{-s/2} \cdot\Vert_H$ for some positive $s$.
}
 \end{remark}
\begin{remark}
{\rm
We point out the essential difference between Theorem~\ref{t1} and
Corollary~\ref{coro02}. This corollary relies on the random squeezing property.
For problems like (\ref{eq-11}) this property is usually proved in the
main space $H$. Therefore the corollary mentioned is applied to functionals
on $H$ only. However in the case of Theorem~\ref{t1} the functionals are
defined on $V$ with $V\subset H$. Thus Theorem~\ref{t1} admits more
singular functionals in comparison with Corollary~\ref{coro02}.
On the other hand Corollary~\ref{coro02} requires convergence of functionals on the
discrete sequence of times $t_n=n$. We note that as in the deterministic case
(see \cite{C4}) it is also possible to consider more general sequences $\{ t_n\}$.
}
\end{remark}

\section{Application to the 2D stochastic Navier-Stokes equations}
We consider the stochastic Navier-Stokes equations
\begin{equation}\label{eq1a}
dv=(-\nu A v+\tilde F(v))dt+dw,
\end{equation}
where
\[
A=-\frac{1}{2}\Delta,\quad \tilde F(v)=-\frac{\nu}{2}\Delta-(v,\nabla)v+f,
\]
as an evolution equation on the rigged Hilbert space $V\subset
H\subset V^\prime$ where $V=\{u\in
\Wsob,\,{\rm
div}\,u=0\}$ and $H=\overline{V}^{L_2(D)}$ is the closure of $V$
in $L^2(D)$. Here $D$ is a bounded domain with sufficiently smooth
boundary  $\partial D$ in $\Bbb{R}^2$, $f\in H$ and  $\nu>0$ is a
constant.  We supplement the Navier-Stokes equations with no-slip
or zero Dirichlet boundary condition $v|_{\partial D} =0$. We
suppose  $w$ is a Wiener process in the space $H^2$ with
covariance $Q$ such that ${\rm tr}_{H^{2}}Q<\infty$. Here and
below we denote by $H^s$ the domain of the operator $A^{s/2}$,
$s>0$. We obviously have $V=H^{1/2}$. In the space $V$ we will use
the norm
$\|\cdot\|_V:=\|\nabla\cdot\|_H=\sqrt{2}\|A^{1/2}\cdot\|_H$.
\\
For different ideas to  treat  this problem, one can find in
\cite{CraFla94},
\cite{CrDeFl97}, \cite{FurVis88}, \cite{BenTem73}, \cite{Vio76}.
\\
We now transform this stochastic equation to a random equation as
in (\ref{eq-11}).
To do this we need a stationary Ornstein~-~Uhlenbeck process $z$.
This process will be generated by the stochastic differential
equation
\begin{equation}\label{eq1a1}
dz+2(k+1)\nu Az\,dt=dw
\end{equation}
for a positive sufficiently large  constant $k$.
It is known (see, e.g. Da Prato and Zabczyk \cite{DaPZab92} Chapter 5) that
there exists a tempered  random variable $z$ in $V$ such that
\[
\Bbb{R}\ni t\to z(\theta_t\omega)
\]
solves (\ref{eq1a1}). This Ornstein~-~Uhlenbeck process $z(\theta_t\omega)$ has
 trajectories in  the space $L^2_{loc}(\Bbb{R};H^{3})$.
The constant $k$ may be considered as a control parameter.
We now consider the nonautonomous differential equation
\begin{equation}\label{eq1b}
\frac{du}{d\,t}+\nu A\, u= F(u,\theta_t\omega),\quad u(0)=x\in H,
\end{equation}
where
\[
F(u,\omega)=-\nu Au-(u\cdot\nabla)u-(z(\omega)\cdot\nabla)u
-(u\cdot\nabla)z(\omega)
\]
\[
-(z(\omega)\cdot\nabla)z(\omega)
+2\nu k Az(\omega) +f.
\]
The idea of this transformation can be found in Crauel and
Flandoli \cite{CraFla94}. Since the coefficients of this equation
have similar properties as the coefficients of the original 2D
Navier~-~ Stokes equations, this equation has a unique solution.
More precisely, we have
\begin{lemma}\label{l3.1}
The solution of (\ref{eq1b}) defines a
continuous
random dynamical system  $\phi$ with respect to $\theta$ on $H$.
Let
\begin{equation}\label{eq1c}
x\to T(\omega,x):=x-z(\omega)
\end{equation}
be a random  homeomorphism on
$H$. Then $T^{-1}(\theta_t\omega,\phi(t,\omega,T(\omega,x)))=:\tilde\phi(t,\omega,x)$
defines a random dynamical system with respect to $\theta$. In particular,
\[
t\to \tilde\phi(t,\omega,x)
\]
solves (\ref{eq1a}).
\end{lemma}
Since $z(\omega)\in V$ the mapping $T$ can be considered as a homeomorphism
on $V$.\\
It is well known  that the random dynamical system $\phi$
has a random compact absorbing forward invariant set $B$ in $V$, see for instance
Crauel and Flandoli \cite{CraFla94}. We now formulate a version of
these results and will prove some additional properties.
\begin{lemma}\label{l3.2}
The random dynamical system $\phi$ has a compact forward invariant
absorbing set $B$ in $H$.
This absorbing set is contained in the closed ball in $H$ with center
zero and with square radius
\begin{equation}\label{eqi}
R^2(\omega)=(1+\eps)\int_{-\infty}^0 m(\theta_\tau\omega)
e^{\nu\lambda_1\tau+\frac{8}{\nu}\int_\tau^0\|z(\theta_s\omega)\|_V^2ds}d\tau,
\end{equation}
where  $\eps>0$ is arbitrary, $\lambda_1$ is the first eigenvalue of the operator
$-\Delta$ with the Dirichlet boundary condition,
\begin{equation}\label{eqi0}
m(\om)=
\frac{4}{\nu}
\left(\frac{2}{\lambda_1}\|z(\omega)\|_V^4+k^2\nu^2\|z(\omega)\|_V^2+
\|f\|_{V^\prime}^2\right)
\end{equation}
and the parameter $k$ in (\ref{eq1a1}) is chosen such that
\begin{equation}\label{eqi1}
\lambda_1>\frac{4\,{\rm tr}_HQ}{(k+1)\nu^3}.
\end{equation}
In addition, $B$ is tempered and
$
t\to\sup_{x\in B(\theta_t\omega)}\|x\|_H^2
$
is a locally integrable stationary process.
\end{lemma}
{\it Proof.$\;$}
We sketch the proof of this lemma.
We obtain by Temam \cite{Tem79} Lemma III.3.4
\[
2\langle  F(u,\theta_t\omega),u\rangle\le
\frac{8}{\nu}\|u\|_H^2\|z(\theta_t\omega)\|_V^2+\frac{8}{\nu\lambda_1}\|z(\theta_t\omega)\|_V^4+
\frac{4}{\nu}\|f\|_{V^\prime}^2+4k^2\nu\|z(\theta_t\omega)\|_V^2.
\]
Let $R_0^2$  be the stationary solution of the random
affine one-dimensional differential equation
\begin{equation}\label{eqm}
\frac{d\rho}{dt}+\nu\lambda_1\rho=
\frac{8}{\nu}\rho\|z(\theta_t\omega)\|_V^2+\frac{8}{\nu\lambda_1}\|z(\theta_t\omega)\|_V^4+
\frac{4}{\nu}\|f\|_{V^\prime}^2+4k^2\nu\|z(\theta_t\omega)\|_V^2.
\end{equation}
This stationary solution $R^2_0(\om)$ exists and it is exponentially attracting
provided
\[
\frac{8}{\nu}\lim_{\tau\to -\infty}\frac{1}{|\tau |}
\int_\tau^0\|z(\theta_s\omega)\|_V^2ds\equiv \frac{8}{\nu}
\Bbb{E}\|z\|_V^2< \nu\lambda_1.
\]
A Simple calculation shows that this relation is equivalent to
(\ref{eqi1}). Moreover  $R^2(\om)=(1+\eps)R^2_0(\om)$  has the
form (\ref{eqi}). The temperedness of $R^2$ follows
 from
  Flandoli and Schmalfu{\ss} \cite{FlaSchm95b} Lemma 7.2. Since the solution
$R^2_0(\theta_t\om)$ of the above equation is continuous,  the mapping
\[
t\to R^2(\theta_t\omega)
\]
is locally integrable. In addition, a comparison argument yields that
the random ball $B(0,R(\omega))$ is forward invariant and forward absorbing.
Finally, we note that
\begin{equation}\label{eq2c}
B(\omega):=\overline{\phi(1,\theta_{-1}\omega,B(0,R(\theta_{-1}\omega)))}\subset B(0,R(\omega))
\end{equation}
is a compact forward invariant and forward absorbing set by the regularization property of $\phi$.
\hfill $\Box$\\

However, there are other compact absorbing sets defined by a ball $B(0,R(\omega))$
with random radius $R(\omega)$, see for instance Flandoli and Langa  \cite{FL}.
In the following we propose another method to calculate moments of (\ref{eqi}).
This technique is based on the standard density of the Girsanov theory.
\begin{lemma}\label{l2}
Let $R^2$ be defined by (\ref{eqi}) then if we choose a $k$ such that
\begin{equation}\label{eq2c1}
\lambda_1>\frac{16\,{\rm tr}_HQ}{(k+1)\nu^3},\qquad
\lambda_1\ge\frac{256\,{\rm tr}_HQ}{(k+1)^2\nu^3},
\end{equation}
we have $\Bbb{E}\,R^8<\infty$.
\end{lemma}
{\it Proof.$\;$}
We rewrite
\[
R^2=(1+\eps)\int_{-\infty}^0
m(\theta_\tau\omega)e^{\nu\lambda_1\tau+c\int_\tau^0\|z\|_V^2}d\tau
\]
for $c=\frac{8}{\nu}$ and some $\eps>0$. We obtain by the Cauchy~-~Schwarz inequality
for an appropriate $c_1>0$
\begin{eqnarray*}
\Bbb{E}R^8&\le& c_1(\Bbb{E}m^8)^\frac{1}{2}
\left(\Bbb{E}\int_{-\infty}^0e^{4\nu\lambda_1\tau
+8c\int_\tau^0\|z\|_V^2}d\tau\right)^\frac{1}{2}
\\
&=& c_1(\Bbb{E}m^8)^\frac{1}{2}
\left(\int^{\infty}_0e^{-4\nu\lambda_1\tau}\cdot
\Bbb{E}e^{8c\int^\tau_0\|z\|_V^2}d\tau\right)^\frac{1}{2}.
\end{eqnarray*}
The first factor is finite, since $z$ is a Gaussian random
variable. We now restrict ourselves to calculate $\Bbb{E}\exp\{
8c\int^\tau_0\|z\|_V^2\}$. Ito's formula applied to
$\|\cdot\|_H^2$ for $z(\theta_t\omega)$ yields:
\[
\|z(\theta_\tau\omega)\|_H^2+2(k+1)\nu\int_0^\tau\|z(\theta_s\omega)\|_V^2ds=\|z(\omega)\|_H^2+
2\int_0^\tau(z,dw)_H+\tau\,{\rm tr}_HQ.
\]
Hence we can derive that
\[
e^{16c\int_0^\tau\|z\|_V^2}\le e^{\frac{8c}{(k+1)\nu}\|z\|_H^2}\cdot
e^{\frac{8c}{(k+1)\nu}\tau\,{\rm tr}_HQ} \cdot e(\tau)\cdot
e^{\frac{256c^2}{(k+1)^2\nu^2}\int_0^\tau(Qz,z)}\, ,
\]
where
\[
e(\tau)\equiv e(\tau,\om)=
\exp\left\{\frac{16c}{(k+1)\nu}\int_0^\tau(z,dw)_H-
\frac{256c^2}{(k+1)^2\nu^2}\int_0^\tau(Qz,z)\right\}\ .
\]
From (\ref{eq2c1}) we have that
$\frac{256c^2{\rm tr}_HQ}{(k+1)^2\nu^2\lambda_1}\le 8c$. Therefore
using the Cauchy~-~ Schwarz inequality
and
\[
(Qz,z)\le{\rm tr}_HQ\|z\|_H^2\le \frac{{\rm tr}_HQ}{\lambda_1}\|z\|_V^2
\]
we have
\begin{equation}\label{eq2c2}
\Bbb{E}e^{8c\int_0^\tau\|z\|_V^2}\le
\left(\Bbb{E}e^{\frac{16c}{(k+1)\nu}\|z\|_H^2}\right)^\frac{1}{2}
\left(\Bbb{E}e(\tau)^2\right)^\frac{1}{2}
e^{\frac{8c}{(k+1)\nu}\tau\,{\rm tr}_HQ}.
\end{equation}
We can use the standard arguments (see, e.g., \cite{Koz78} and
\cite{GihSko79}) to find the estimate $\Bbb{E}[e(\tau)^2]\le 1$
for the mean value of Girsanov's density $e(\tau)^2$. The value
 $z(\omega)$ is a Gaussian variable in $H$ with the zero mean
and with the covariance
\[
\Bbb{E}\langle z, h_1\rangle\langle z, h_2\rangle =\langle \tilde Qh_1,
h_2\rangle,\quad h_1, \, h_2\in H,
\]
where
\[
\tilde{Q}=\int_0^\infty e^{-2t(k+1)\nu A}Qe^{-2t(k+1)\nu A}\, dt.
\]
Therefore simple calculation (see, e.g. Kuo~\cite{K}) Page 105
shows that the first factor in the right hand side of
(\ref{eq2c2}) is finite provided $\frac{16c}{(k+1)\nu}< \frac{1}{2
{\rm tr}_H\tilde{Q}}$. Moreover
\[
\Bbb{E}e^{\frac{16c}{(k+1)\nu}\|z\|_H^2}\le \exp\left\{
\frac{16c{\rm tr}_H\tilde Q}{(k+1)\nu -32c{\rm tr}_H\tilde
Q}\right\}\ .
\]
However it is easy to see that
${\rm tr}_H\tilde{Q}\le \frac{{\rm tr}_HQ}{2\lambda_1(k+1)\nu}$.
Therefore from the second assumption of (\ref{eq2c1}) we have
that
\begin{equation}\label{eq23c}
\Bbb{E}e^{\frac{16c}{(k+1)\nu}\|z\|_H^2}\le \exp\left\{
\frac{8c{\rm tr}_H Q}{\lambda_1(k+1)^2\nu^2 -16c{\rm tr}_H
Q}\right\} <\infty.
\end{equation}
Since from (\ref{eq2c1}) we also have
$4\nu\lambda_1>\frac{8c{\rm tr}_HQ}{(k+1)\nu}$, the expectation of $R^8$ is
finite.
\hfill $\Box$\\

The following lemma allows us to conclude the existence of
a set ${\cal L}$ of determining functionals for
the random dynamical system $\tilde\phi$ generated by (\ref{eq1a}) if
the random dynamical system $\phi$ generated by (\ref{eq1b})
has the set of determining functionals ${\cal L}$. This
lemma is formulated for more general transformations than
(\ref{eq1c}).
\begin{lemma}
Suppose that the random dynamical systems $\tilde\phi$ and $\phi$
are conjugated by a random homeomorphism $T$ on $H$, i.e.
$\tilde\phi(t,\om,\tilde x(\om))=
T^{-1}(\theta_t\om,\phi(t,\om, x(\om))$, where $x(\om)=T(\om,\tilde x(\om))$.
Suppose that $\phi$ has a compact absorbing set and forward invariant
set $B$. Then
$\tilde\phi(t,\omega,\tilde x_1(\omega))-\tilde\phi(t,\omega,\tilde x_2(\omega))$
tends to zero in probability for $t\to\infty$ if and only if
$\phi(t,\omega,x_1(\omega))-\phi(t,\omega,x_2(\omega))$
tends to zero in probability for $t\to\infty$.
Here $x_i(\om)=T(\om,\tilde x_i(\om))$.
\end{lemma}
{\it Proof.$\;$} Suppose that
$\phi(t,\omega,x_1(\omega))-\phi(t,\omega,x_2(\omega))$ tends to
zero in probability for $t\to\infty$. By the absorbing property of
$B$ we can assume that $x_1(\omega),\,x_2(\omega)\in B(\omega)$.
For any $\eps>0$ there exists a compact set $C_\eps$ such that
$C_\eps\supset B(\omega)$ with probability bigger than $1-\eps$.
Indeed, this follows by the regularization property of $\phi$ and
by the construction of $B$ in (\ref{eq2c}). $T^{-1}(\omega)$ is
uniformly continuous on $C_\eps$: for any $\omega\in\Omega$,
$\mu>0$, $y_1,\,y_2\in C_\eps$ there exists a $\delta(\omega)>0$
such that if $\|y_1-y_2\|_H<\delta(\omega)$ then
$\|T^{-1}(\omega,y_1)-T^{-1}(\omega,y_2)\|_H<\mu$. On the other
hand since $\delta(\omega)>0$ there exists a $\delta_\eps>0$:
\[
\Bbb{P}(\delta_\eps<\delta(\omega))>1-\eps.
\]
Hence for
sufficiently large $t_\eps$ we have
\begin{eqnarray*}
&&
\Bbb{P}(\|\phi(t,\omega,x_1(\omega))-\phi(t,\omega,x_2(\omega))\|_H
>\delta(\theta_t\omega))\\
&& \qquad \le \eps+\Bbb{P}(\|\phi(t,\omega,x_1(\omega))
-\phi(t,\omega,x_2(\omega))\|_H
>\delta_\eps)<2\eps
\end{eqnarray*}
if $t\ge t_\eps$.
Hence
\begin{eqnarray*}
&&
\|\tilde\phi(t,\omega,\tilde x_1(\omega))-\tilde\phi(t,\omega,\tilde x_2(\omega))\|_H\\
&&\qquad=
\|T^{-1}(\theta_t\omega,\phi(t,\omega,x_1(\omega)))-
T^{-1}(\theta_t\omega,\phi(t,\omega,x_2(\omega)))\|_H<\mu
\end{eqnarray*}
with probability bigger than $1-3\eps$ for  $t>t_\eps$.\\
$\tilde B=T^{-1}(B)$ is a compact forward invariant absorbing set
for (\ref{eq1a}) if and only if $B$ is a compact forward invariant absorbing set
for (\ref{eq1b}). Therefore we can show the second direction similarly as the
proof above for  the first direction.
\hfill $\Box$\\

\begin{coro}\label{coro1}
The set of linear functionals  ${\cal L}$ on $V$ is determining
in probability for the random dynamical system $\tilde\phi$
generated by (\ref{eq1a}) if and only if ${\cal L}$ is determining
in probability for $\phi$ defined by (\ref{eq1b}).
\end{coro}
{\it Proof.$\;$} The proof is based on the fact that for some
$l\in{\cal L}$ the limit in probability for $t\to\infty$ of
$l(\tilde\phi(t,\omega,\tilde
x_1(\omega))-\tilde\phi(t,\omega,\tilde x_2(\omega)))$ is zero if
and only if
$l(\phi(t,\omega,x_1(\omega))-\phi(t,\omega,x_2(\omega)))$ tends
to zero in probability for $t\to\infty$ which follows from the
particular shape of $T$. On the other hand, we can
 also
apply the last lemma.
\hfill $\Box$\\

For the following we need two a priori estimates for $\phi$:
\begin{lemma}\label{l3}
The random dynamical system $\phi$ satisfies the
following a priori estimate
\begin{eqnarray*}
&&\nu\sup_{x\in B(\omega)} \int_0^t\|\phi(\tau,\omega,
x)\|_V^2d\tau
\le R^2(\omega)+\frac{8}{\nu}\int_0^t\|z(\theta_\tau\omega)\|_H^2\|z(\theta_\tau\omega)\|_V^2d\tau\\
&&\quad +\frac{4}{\nu}t\|f\|_{V^\prime}^2+\nu
k^2\int_0^t\|z(\theta_\tau\omega)\|_V^2d\tau
+c_E M\int_0^t\|z(\theta_\tau\omega)\|_{H^3}^2d\tau\\
&&\quad+\frac{c_E}{M}\int_0^t R^4(\theta_\tau\omega)d\tau.
\end{eqnarray*}
where $M$ is an arbitrarily positive number and $c_E$ is the norm of the
 embedding operator of $H^3=D(A^{3/2})$ into  the space $W^1_\infty (D)$
of two-dimensional functions $v$ such that $v, \nabla v\in L^\infty (D)$.
Similarly, for an appropriate polynomial $p$
\begin{eqnarray*}
&&2\nu\sup_{x\in  B(\omega)}
\int_0^t\|\phi(\tau,\omega,x)\|_H^2\|\phi(\tau,\omega,x)\|_V^2d\tau
\le R^4(\omega)\\
&&\qquad +\int_0^t
p(\|z(\theta_\tau\omega)\|_H^2,\|z(\theta_\tau\omega)\|_V^2,
\|z(\theta_\tau\omega)\|_H^3, \|f\|_{V^\prime}^2)d\tau
+\int_0^tR^8(\theta_\tau\omega)d\tau \ .
\end{eqnarray*}
\end{lemma}
This a priori estimate is based on the calculation of $\|u(t)\|_H^2$ for
(\ref{eq1b}).
The   term
$\langle(u\cdot\nabla) z),u\rangle$ arising in the  calculation  can be estimated by the Sobolev lemma:
\[
|\langle (u\cdot\nabla) z,u\rangle|\le c_E\|z\|_{H^3}\|u\|_H^2\le \frac{c_EM}{2}\|z\|_{H^3}^2
+\frac{c_E}{2M}R^4
\]
because $x\in B$. The second estimate follows similarly for $\|u(t)\|_H^4$.\\

\begin{lemma}\label{l4} Under conditions (\ref{eq2c1})
the following estimate holds:
\begin{eqnarray}\label{eqn01}
\Sigma_k&\equiv\lim\sup_{m\to\infty}\frac{1}{m}{\Bbb
E}\left\{
 \sup_{x\in B(\omega)}
\int_0^m\|\phi(\tau,\omega, x)\|_V^2d\tau\right\}\\
&\le \left(\frac{4}{\nu^2}\|f\|_{V^\prime}^2+
g_k(\nu,\lambda_1, Q)\right)\cdot
\left( 1+
h_k(\nu,A, Q)\right),\nonumber
\end{eqnarray}
where
\begin{equation}\label{eqn01a}
g_k(\nu,\lambda_1, Q)=
a_0 \frac{{\rm tr}_HQ}{\nu}\cdot\left( k+  \frac{a_1{\rm tr}_HQ}
{(k+1)^2\lambda_1\nu^3}\right)
\end{equation}
and
\begin{equation}\label{eqn01b}
h_k(\nu,A, Q)= 2c_E \left( \frac{ [{\rm tr}_HQA^2]^2}{
\nu^3\lambda_1^3(k+1)\cdot [\nu^3\lambda_1(k+1) -16 {{\rm tr}_H
Q}]}\right)^{1/4}.
\end{equation}
Here $a_0$ and $a_1$ are some absolute constants and $c_E$ is the same as in
Lemma~\ref{l3}.
\end{lemma}
{\it Proof.$\;$}
It follows from Lemma~\ref{l3} that
\[
\Sigma_k\le \frac{4}{\nu^2}\|f\|_{V^\prime}^2+
\frac{8}{\nu^2}{\Bbb E}\left( \|z\|_H^2\|z\|_V^2\right)+ k^2 {\Bbb
E}\|z\|_V^2 +\frac{c_E}{\nu} M {\Bbb E}\|z\|_{H^3}^2
+\frac{c_E}{\nu M}{\Bbb E} R^4.
\]
If we choose $M=\left( {\Bbb E}\|z\|_{H^3}^2\right)^{-1/2}\cdot
\left( {\Bbb E} R^4\right)^{1/2}$, then we obtain
\begin{eqnarray*}
\Sigma_k&\le \frac{4}{\nu^2}\|f\|_{V^\prime}^2+
\frac{8}{\nu^2}\left({\Bbb E}
\|z\|_H^4\right)^{1/2}\cdot\left( {\Bbb E}\|z\|_V^4\right)^{1/2}\\
& + k^2 {\Bbb E}\|z\|_V^2 +\frac{2c_E}{\nu} \left( {\Bbb
E}\|z\|_{H^3}^2\right)^{1/2} \cdot\left( {\Bbb E}
R^4\right)^{1/2}.
\end{eqnarray*}
Using the definition of $z$ it is easy to find that for any positive
$\alpha\in [0,3/2]$ we have
\begin{equation}\label{eqn01c}
{\Bbb E}\|A^{\alpha}z\|_{H}^2=\frac{1}{4(k+1)\nu}{\rm
tr}_H(QA^{2\alpha-1}).
\end{equation}
Furthermore  it is clear that
\begin{equation}\label{eqn02}
{\Bbb E}\|A^{\alpha}z\|_{H}^{2l}\le c_l \left( {\Bbb
E}\|A^{\alpha}z\|_{H}^2\right)^l,\quad l=1,2,\ldots,
\end{equation}
with appropriate constants $c_l$. Therefore we have
\[
\Sigma_k\le \frac{4}{\nu^2}\|f\|_{V^\prime}^2+ \frac{c_0[{\rm
tr}_HQ]^2}{\nu^4 (k+1)^2\lambda_1}+ \frac{k}{2\nu}\cdot {\rm
tr}_HQ+ \frac{c_E[{\rm tr}_HQA^2]^{1/2}}{\nu^{3/2} (k+1)^{1/2}}
\cdot\left( {\Bbb E} R^4\right)^{1/2}.
\]
with some absolute constant $c_0$. Now we estimate $\left( {\Bbb
E} R^4\right)^{1/2}$. We use the idea of the proof of
Lemma~\ref{l2}. It is clear that
\[
\left( {\Bbb E} R^4\right)^{1/2}\le
(1+\eps)\left(\frac{3}{2\nu\lambda_1}\right)^{3/4} \left( {\Bbb E}
m^4\right)^{1/4} \left(
  {\Bbb E}\int_{-\infty}^0
e^{2\nu\lambda_1\tau+4c\int_\tau^0\|z(\theta_s\omega)\|_V^2ds}d\tau
\right)^{1/4},
\]
where $m(\om)$ is given by (\ref{eqi0})  and $c=8/\nu$.
Using  Girsanov's trick and (\ref{eq2c2}) and (\ref{eq23c}) with $c:=c/2$
we have
\[
\Bbb{E}e^{4c\int_0^\tau\|z\|_V^2}\le \exp\left\{ \frac{2c{\rm
tr}_H Q}{\lambda_1(k+1)^2\nu^2 -8c{\rm tr}_H Q}\right\} \cdot
\exp\left\{ \frac{4c}{(k+1)\nu}\tau\,{\rm tr}_HQ\right\}
\]
under conditions (\ref{eq2c1}). However (\ref{eq2c1}) implies
that $\lambda_1(k+1)^2\nu^2 \ge 32c{\rm tr}_H Q$. Therefore
\begin{equation}\label{eqn03}
\Bbb{E}e^{4c\int_0^\tau\|z\|_V^2}\le \exp\left\{ \frac{1}{12}+
\frac{4c}{(k+1)\nu}\tau\,{\rm tr}_HQ\right\}
\end{equation}
under conditions (\ref{eq2c1}). From  (\ref{eqn03}) we have
\[
\left( {\Bbb E} R^4\right)^{1/2}\le (1+\eps)
e^{1/48}\left(\frac{3}{2\nu\lambda_1}\right)^{3/4} \cdot
\left(2\nu\lambda_1-\frac{32}{(k+1)\nu^2} {\rm tr}_H
Q\right)^{-1/4} \left( {\Bbb E} m^4 \right)^{1/4}.
\]
Now we estimate $\left( {\Bbb E} m^4 \right)^{1/4}$. It is clear
from (\ref{eqi0}) that
\[
\left( {\Bbb E} m^4 \right)^{1/4}\le \frac{4}{\nu}
\left\{\frac{2}{\lambda_1} \left( {\Bbb
E}\|z(\omega)\|_V^{16}\right)^{1/4}+ k^2\nu^2\left(  {\Bbb E}
\|z(\omega)\|_V^8\right)^{1/4}+ \|f\|_{V^\prime}^2\right\}.
\]
Therefore using (\ref{eqn02}) we obtain
\[
\left( {\Bbb E} m^4 \right)^{1/4}\le \frac{4}{\nu}
\left\{\|f\|_{V^\prime}^2+ \frac{c_1[{\rm tr}_HQ]^2}{\nu^2
(k+1)^2\lambda_1}+ c_2k\nu\cdot {\rm tr}_HQ\right\},
\]
where $c_1$ and $c_2$ are absolute constants. Put   all these estimates together
we obtain the upper bound (\ref{eqn01}) for $\Sigma_k$.
\hfill $\Box$\\

We have seen that $B$, defined in (\ref{eq2c}), is bounded in $V$, and hence it is a compact set in $H$
which is tempered with respect the $H$ norm.
We now prove that $B$ is also tempered and locally integrable in $V$.
\begin{lemma}
The random variable $\sup_{x\in  B(\omega)}\|x\|_V^2$ is tempered
and the mapping\\ $ t\to\sup_{x\in  B(\theta_t\omega)}\|x\|_V^2 $
is locally integrable.
\end{lemma}
{\it Proof.$\;$}
To obtain an estimate in $V$ we use the standard method which is based on
the formula
\[
\frac{d}{dt}(t\|u(t)\|_V^2)=\|u(t)\|_V^2+t\frac{d\|u(t)\|_V^2}{dt}.
\]
for $t=1$ by Temam \cite{Tem79} Lemma III.3.8 and
\[
|((u\cdot\nabla) v,w)_V|\le c \|u\|_H^{\frac{1}{2}}\|u\|_V^{\frac{1}{2}}
\|v\|_V^{\frac{1}{2}}\|Av\|_H^{\frac{1}{2}}\|Aw\|_H
\]
for sufficiently regular $u,\,v,\,w$ and $c>0$ which allows us to write
by (\ref{eq1b})
\begin{eqnarray*}
\frac{d}{dt}(t\|u(t)\|_V^2)&\le&K(\|u(t)\|_H^2\|u(t)\|_V^2
+\|u(t)\|_H^2+
\|z(\theta_t\omega)\|_H^2\|z(\theta_t\omega)\|_V^2) \\
&\times &(t\|u(t)\|_V^2)\\
&+&p(\|f\|_{V^\prime},\|z(\theta_t\omega)\|_H,\|z(\theta_t\omega)\|_V,\|Az(\theta_t\omega)\|_H)+\|u(t)\|_V^2
\end{eqnarray*}
where $p$ is an appropriate polynomial and $K$ an appropriate positive
constant. Consequently, by the Gronwall lemma
\begin{eqnarray*}
\sup_{x\in B}\|\phi(1,\omega,x)\|_V^2&\le&
\exp\left(K\sup_{x\in B}\int_0^1\|\phi(\tau,\omega,x)\|_H^2d\tau\right) \\
&\times&\exp\left(K\sup_{x\in  B}\int_0^1
\|\phi(\tau,\omega,x)\|_H^2\|\phi(\tau,\omega,x)\|_V^2d\tau\right) \\
&\times&\exp\left(K\int_0^1\|z(\theta_\tau\omega)\|_V^2
\|z(\theta_\tau\omega)\|_H^2d\tau\right) \\
&\times&\left(\int_0^1p(\tau)d\tau+\sup_{x\in
B}\int_0^1\|\phi(\tau,\omega,x)\|_V^2d\tau\right).
\end{eqnarray*}
Note that a product of random variables is tempered if each factor
is tempered. To see that the first factor
of the right hand side is tempered we use the estimate
\[
\Bbb{E}\sup_{\scriptsize
\begin{array}{c}
x\in  B\\
s\in[0,1]
\end{array}
}\int_0^1\|\phi(\tau,\theta_s\omega,x)\|_H^2d\tau
\le\Bbb{E}\int_0^2R^2(\theta_s\omega)ds= 2\Bbb{E}R^2<\infty
\]
by Lemma \ref{l2} and the forward invariance of $B$, see Arnold \cite{Arn98} Proposition 4.1.3.
Similarly, we get for the next factor
\[
\Bbb{E}\sup_{\scriptsize
\begin{array}{c}
x\in  B\\
s\in[0,1]
\end{array}
}
\int_0^1\|\phi(\tau,\theta_s\omega,x)\|_H^2\|\phi(\tau,\theta_s\omega,x)\|_V^2d\tau<\infty
\]
which follows from Lemma \ref{l3}. However, to justify this
estimate we also need that $\Bbb{E}\sup_{s\in
[0,1]}R^4(\theta_s\omega)<\infty$. For this expression we obtain
an estimate if we calculate in (\ref{eqm})
$R_0^4(\theta_s\omega)=\rho^2(s)$ by the chain rule. Then we can
estimate this supremum by $R^4(\omega)$ and some integrals of
norms from $z$ which have a finite expectation.
The temperedness of the remaining factors follow similarly.\\
The local integrability follows by the continuity of $t\to R(\theta_t\omega)$
and the local integrability of the norms of $z$.
\hfill $\Box$\\

We   are now in a position to formulate the main theorem of this section.
\begin{theorem}\label{t2}
Let ${\cal L}$ be a set of linear functionals on $V$ with
completeness defect $\eps_{\cal L}$. Assume that for some $k$
satisfying  (\ref{eq2c1}) the  completeness defect
$\eps_{\cal L}$ possesses the property
\begin{equation}\label{eqn1}
\frac{4}{\nu}\cdot \left(\frac{4}{\nu^2}\|f\|_{V^\prime}^2+
g_k(\nu,\lambda_1, Q)\right)\cdot
\left( 1+
h_k(\nu, A,Q)\right)+
\frac{2}{(k+1)\nu}{\rm tr}_HQ
< \nu\eps_{\cal L}^{-2},
\end{equation}
where $g_k(\nu,\lambda_1, Q)$ and $h_k(\nu,A, Q)$ are given by (\ref{eqn01a})
and  (\ref{eqn01b}).
Then ${\cal L}$ is a system of determining functionals
in probability for the 2D stochastic Navier-Stokes equation  (\ref{eq1a}).
\end{theorem}
{\it Proof.$\;$}
We are going to apply Theorem \ref{t1}. The temperedness and local
integrability of $\sup_{x\in B(\theta_t\om}\| x\|^2_V$ follow by
the last lemma. Then we get the  assertion if we choose $m$
sufficiently large. Indeed, in the case of large $m$ we can reduce
the influence of $\Bbb{E} R^2$. By Corollary \ref{coro1} it is
sufficient to show that ${\cal L}$ is a set of determining
functionals for $\phi$ generated by (\ref{eq1b}). The properties
of $F$ allow us to estimate the Lipschitz constant
\[
l(x_1,x_2,\omega)=
\frac{2}{\nu}(\|x_1\|_V^2+\|z(\omega)\|_V^2).
\]
The measurability of $l$ follows straightforwardly.
We should also take $c=\frac{\nu}{2}$ in (\ref{eqb1}). Therefore we can apply
Theorem~\ref{t1} if
\begin{equation}\label{eqn2}
\frac{4}{\nu m}{\Bbb E}\left\{
 \sup_{x\in B(\omega)}
\int_0^m\|\phi(\tau,\omega, x)\|_V^2d\tau\right\}+
\frac{4}{\nu}{\Bbb E}\| z\|^2_V<\nu\eps_{\cal L}^{-2}
\end{equation}
for some $m$. We can find $m$ with the property (\ref{eqn2}), if
\[
\frac{4}{\nu}\Sigma_k+ \frac{4}{\nu}{\Bbb E}\|
z\|^2_V<\nu\eps_{\cal L}^{-2}.
\]
The last relation follows from Lemma~\ref{l4}, the relation (\ref{eqn01c})
and  (\ref{eqn1}).
\hfill $\Box$\\

\begin{remark}
{\rm
In the limit ${\rm tr}_{H}QA^2\to 0$
relation (\ref{eqn1}) turns in the inequality
\begin{equation}\label{eqn3}
\eps_{\cal L}<\frac{4\nu^2}{\|f\|_{V^\prime}}.
\end{equation}
Thus if the estimate (\ref{eqn3}) is valid, then there exists a  constant
$\delta_0>0$ such that under condition ${\rm tr}_{H}QA^2<\delta_0$
the set  ${\cal L}$ is a set of determining functionals
in probability for the 2D stochastic Navier-Stokes equations.
We also note that  estimate (\ref{eqn3}) is the same order as the best
known estimate for the completeness defect in the case of deterministic
$2D$ Navier~-~Stokes equations with the periodic boundary conditions
(see the survey \cite{C3} and the references therein).
However in the last case relation (\ref{eqn3}) involves the completeness defect
with respect to the pair $D(A)$ and $H$ and it leads to better estimates for
the number of determining functionals.
}
\end{remark}
\begin{remark}
{\rm As an application of Theorem~\ref{t1} and \ref{t1a} we can
consider the equation
\[
\partial_t u =\Delta u -f(u) +\partial_t W(t,\om)
\]
in a bounded domain, where $f(u)$ is a polynomial of odd degree
with positive leading coefficient. We can also consider 2D
stochastic Navier-Stokes equations with multiplicative white noise
$u\,dw$, where $w$ is a scalar Wiener process. In this case  we
have to use the transformation $T(\omega,x)=x\,e^{-z(\omega)}$
where $z$ defines a one dimensional stationary Ornstein-Uhlenbeck
process generated by $dz+z\,dt=dw$. This equation has been
investigated for instance in Schmalfu{\ss} \cite{Schm95a} but with
a little bit different transformation $T$. }

\end{remark}

\medskip

{\bf Acknowledgment.} A part of this work was done at the
Oberwolfach Mathematical Research Institute, Germany, while J.
Duan and B. Schmalfu{\ss} were Research in Pairs Fellows,
supported by the {\em Volkswagen Stiftung}.

\end{document}